\documentstyle[amssymb,11pt]{article}
\title{\bf Symmetries of Surface Singularities}
\author{Gerd M\"uller}
\date{March 1996}
\begin{document}
\newcommand{\sub}{\subseteq}
\newcommand{\Aut}{{\rm Aut}}
\newcommand{\C}{{\Bbb C}}
\newcommand{\Z}{{\Bbb Z}}
\renewcommand{\P}{{\Bbb P}}
\newcommand{\GL}{{\rm GL}}
\newcommand{\PSL}{{\rm PSL}}
\newcommand{\im}{{\rm im}}
\renewcommand{\O}{{\cal O}}
\newcommand{\G}{\Gamma}
\renewcommand{\a}{\alpha}
\renewcommand{\b}{\beta}
\newcommand{\g}{\gamma}
\renewcommand{\l}{\lambda}
\newcommand{\s}{\sigma}
\renewcommand{\t}{\tau}
\newcommand{\ep}{\varepsilon}
\renewcommand{\L}{\Lambda}
\newcommand{\w}{\omega}
\newcommand{\T}{{\rm T}}

\newcommand{\Aone}[1]{
  \unitlength0.8mm
  \begin{picture}(0,9)
  \put(0,6){\circle*{2}}
  \put(0,0){\makebox(0,0)[b]{$#1$}}
  \end{picture}
  }

\newcommand{\Atwo}[2]{
  \unitlength0.8mm
  \begin{picture}(20,9)
  \multiput(0,6)(20,0){2}{\circle*{2}}
  \put(0,0){\makebox(0,0)[b]{$#1$}}
  \put(20,0){\makebox(0,0)[b]{$#2$}}
  \put(0,6){\line(1,0){20}}
  \end{picture}
  }

\newcommand{\Athree}[3]{
  \unitlength0.8mm
  \begin{picture}(40,9)
  \multiput(0,6)(20,0){3}{\circle*{2}}
  \put(0,0){\makebox(0,0)[b]{$#1$}}
  \put(20,0){\makebox(0,0)[b]{$#2$}}
  \put(40,0){\makebox(0,0)[b]{$#3$}}
  \put(0,6){\line(1,0){40}}
  \end{picture}
  }

\newcommand{\Afour}[4]{
  \unitlength0.8mm
  \begin{picture}(60,9)
  \multiput(0,6)(20,0){4}{\circle*{2}}
  \put(0,0){\makebox(0,0)[b]{$#1$}}
  \put(20,0){\makebox(0,0)[b]{$#2$}}
  \put(40,0){\makebox(0,0)[b]{$#3$}}
  \put(60,0){\makebox(0,0)[b]{$#4$}}
  \put(0,6){\line(1,0){60}}
  \end{picture}
  }

\newcommand{\Afive}[5]{
  \unitlength0.8mm
  \begin{picture}(80,9)
  \multiput(0,6)(20,0){5}{\circle*{2}}
  \put(0,0){\makebox(0,0)[b]{$#1$}}
  \put(20,0){\makebox(0,0)[b]{$#2$}}
  \put(40,0){\makebox(0,0)[b]{$#3$}}
  \put(60,0){\makebox(0,0)[b]{$#4$}}
  \put(80,0){\makebox(0,0)[b]{$#5$}}
  \put(0,6){\line(1,0){80}}
  \end{picture}
  }

\newcommand{\Dfour}[4]{
  \unitlength0.8mm
  \begin{picture}(40,29)
  \multiput(0,6)(20,0){3}{\circle*{2}}
  \put(20,26){\circle*{2}}
  \put(0,0){\makebox(0,0)[b]{$#2$}}
  \put(20,0){\makebox(0,0)[b]{$#3$}}
  \put(40,0){\makebox(0,0)[b]{$#4$}}
  \put(23,26){\makebox(0,0)[l]{$#1$}}
  \put(0,6){\line(1,0){40}}
  \put(20,6){\line(0,1){20}}
  \end{picture}
  }

\newcommand{\Dfive}[5]{
  \unitlength0.8mm
  \begin{picture}(60,29)
  \multiput(0,6)(20,0){4}{\circle*{2}}
  \put(20,26){\circle*{2}}
  \put(0,0){\makebox(0,0)[b]{$#2$}}
  \put(20,0){\makebox(0,0)[b]{$#3$}}
  \put(40,0){\makebox(0,0)[b]{$#4$}}
  \put(60,0){\makebox(0,0)[b]{$#5$}}
  \put(23,26){\makebox(0,0)[l]{$#1$}}
  \put(0,6){\line(1,0){60}}
  \put(20,6){\line(0,1){20}}
  \end{picture}
  }

\newcommand{\Dsix}[6]{
  \unitlength0.8mm
  \begin{picture}(80,29)
  \multiput(0,6)(20,0){5}{\circle*{2}}
  \put(20,26){\circle*{2}}
  \put(0,0){\makebox(0,0)[b]{$#2$}}
  \put(20,0){\makebox(0,0)[b]{$#3$}}
  \put(40,0){\makebox(0,0)[b]{$#4$}}
  \put(60,0){\makebox(0,0)[b]{$#5$}}
  \put(80,0){\makebox(0,0)[b]{$#6$}}
  \put(23,26){\makebox(0,0)[l]{$#1$}}
  \put(0,6){\line(1,0){80}}
  \put(20,6){\line(0,1){20}}
  \end{picture}
  }

\newcommand{\Esix}[6]{
  \unitlength0.8mm
  \begin{picture}(80,29)
  \multiput(0,6)(20,0){5}{\circle*{2}}
  \put(40,26){\circle*{2}}
  \put(0,0){\makebox(0,0)[b]{$#2$}}
  \put(20,0){\makebox(0,0)[b]{$#3$}}
  \put(40,0){\makebox(0,0)[b]{$#4$}}
  \put(60,0){\makebox(0,0)[b]{$#5$}}
  \put(80,0){\makebox(0,0)[b]{$#6$}}
  \put(43,26){\makebox(0,0)[l]{$#1$}}
  \put(0,6){\line(1,0){80}}
  \put(40,6){\line(0,1){20}}
  \end{picture}
  }

\newcommand{\Eseven}[7]{
  \unitlength0.8mm
  \begin{picture}(100,29)
  \multiput(0,6)(20,0){6}{\circle*{2}}
  \put(40,26){\circle*{2}}
  \put(0,0){\makebox(0,0)[b]{$#2$}}
  \put(20,0){\makebox(0,0)[b]{$#3$}}
  \put(40,0){\makebox(0,0)[b]{$#4$}}
  \put(60,0){\makebox(0,0)[b]{$#5$}}
  \put(80,0){\makebox(0,0)[b]{$#6$}}
  \put(100,0){\makebox(0,0)[b]{$#7$}}
  \put(43,26){\makebox(0,0)[l]{$#1$}}
  \put(0,6){\line(1,0){100}}
  \put(40,6){\line(0,1){20}}
  \end{picture}
  }

\maketitle
\section{Introduction}
The study of reductive group actions on a normal surface singularity $X$
is facilitated by the fact that the group $\Aut\: X$ of automorphisms of
$X$ has a maximal reductive algebraic subgroup $G$ which contains every
reductive algebraic subgroup of $\Aut\: X$ up to conjugation. If $X$ is
not weighted homogeneous then this maximal group $G$ is finite (Scheja,
Wiebe). It has been determined for cusp singularities by Wall. On the
other hand, if $X$ is weighted homogeneous but not a cyclic quotient
singularity then the connected component $G_1$ of the unit coincides
with the $\C^*$ defining the weighted homogeneous structure (Scheja,
Wiebe and Wahl). Thus the main interest lies in the finite group
$G/G_1$. Not much is known about $G/G_1$. Ganter has
given a bound on its order valid for Gorenstein singularities which are
not log-canonical. Aumann-K\"orber has determined $G/G_1$ for all
quotient singularities.
\\[1.5ex]
We propose to study $G/G_1$ through the action of $G$ on the minimal good
resolution $\tilde{X}$ of $X$. If $X$ is weighted homogeneous but not
a cyclic quotient singularity let $E_0$ be the central curve of the
exceptional divisor of $\tilde{X}$. We show that the natural homomorphism
$G\to\Aut\: E_0$ has kernel $\C^*$ and finite image. In particular, this
reproves the result of Scheja, Wiebe and Wahl mentioned above.
Moreover, it allows to view $G/G_1$ as a subgroup of
$\Aut\: E_0$. For simple elliptic singularities it equals
$(\Z_b\times\Z_b)\rtimes\Aut_0\, E_0$ where $-b$ is the self intersection
number of $E_0$, $\Z_b\times\Z_b$ is the group of $b$-torsion points of the
elliptic curve $E_0$ acting by translations,
and $\Aut_0\, E_0$ is the group of automorphisms
fixing the zero element of $E_0$. If $E_0$ is rational then $G/G_1$ is the
group of automorphisms of $E_0$ which permute the intersection points
with the branches of the exceptional divisor while preserving the Seifert
invariants of these branches. When there are exactly three branches we
conclude that $G/G_1$ is isomorphic to the group of automorphisms of the
weighted resolution graph. This applies to all non-cyclic quotient
singularities as well as to triangle singularities. We also
investigate whether the maximal reductive automorphism group is a direct
product $G\simeq G_1\times G/G_1$. This is the case, for instance, if
the central curve $E_0$ is rational of even self intersection number or if
$X$ is Gorenstein such that its nowhere zero $2$-form $\omega$ has degree
$\pm 1$. In the latter case there is a ``natural'' section
$G/G_1\hookrightarrow G$ of $G\to G/G_1$ given by the group of automorphisms
in $G$ which fix $\omega$. For a simple elliptic singularity one has
$G\simeq G_1\times G/G_1$ if and only if $-E_0\cdot E_0=1$.
\\[1.5ex]
In the weighted homogeneous case, we show that $G/G_1$ acts
faithfully on the homology $H_1(L,\Z)$ of the link of $X$. For a hypersurface
in $\C^3$, not an $A_k$-singularity, defined by a weighted homogeneous
polynomial $f$ this will be rephrased in terms of the group $H$ of linear
right equivalences of $f$. Namely, with $F$ denoting the Milnor fibre of $f$,
the group $H$ acts faithfully on the Milnor lattice $H_2(F,\Z)$ and
intersects the monodromy group in the cyclic group generated by the
monodromy operator.
\\[1.5ex]
The results of this paper were obtained during a stay at the Mathematical
Institute of Utrecht University. I thank its members for their hospitality
and the DFG for financial support. I profited from discussions with
Theo de Jong, Felix Leinen, Eduard Looijenga, Peter Slodowy, Tonny Springer,
Joseph Steenbrink and Jan Stevens.
Especially I wish to thank Dirk Siersma
for his suggestions and support.
\section{The action on the central curve}
Throughout this paper let $X=(X,0)$ denote a normal surface singularity
with local ring $(\O_X,m_X)$.
We are interested in algebraic subgroups $G$ of the group $\Aut\: X$ of
automorphisms of $X$. This means that $G\le\Aut\: X$ is an abstract
subgroup equipped with the structure of an algebraic group such that
the natural representations of $G$ on all higher cotangent spaces
$m_X^k/m_X^{k+1}$ of $X$ are rational.
There is an abundance of unipotent
algebraic subgroups of $\Aut\: X$, see \cite[Theorem 5]{MMonats}.
But the situation
becomes simpler if one restricts to reductive groups:
\begin{trivlist}
  \item[] \bf Theorem 1. \it
  There is a maximal reductive algebraic subgroup $G\le\Aut\: X$
  containing every
  reductive algebraic subgroup of $\Aut\: X$ up to conjugation.
  \item[] Proof. \rm
  This follows from \cite[Theorem 1]{HM} since a normal surface singularity
  (just as any isolated singularity) can be defined by polynomials, see
  \cite[Theorem 3.8]{A}. If one wants to avoid the use of
  \cite[Theorem 1]{HM}, which depends on a deep result of Popescu and
  Rotthaus, one can argue more elementary as follows: Let $\hat{X}$ be
  the algebroid space corresponding to the formal completion $\widehat{\O_X}$
  of $\O_X$.
  By \cite[Satz 4]{MCrelle}
  the group of automorphisms of $\hat{X}$
  contains a maximal reductive subgroup $G$. It is enough to make $G$ act
  on $X$ itself, with the same representation on the cotangent space,
  \cite[Lemma]{HM}. If $G$ is finite then \cite[Remark on p.\ 184]{HM}
  applies. Otherwise, $G$ contains some $\C^*$. It follows from a result of
  Scheja and Wiebe \cite[3.1]{SW} that in this case the normal surface
  $X$ in fact admits a good $\C^*$-action, i.\ e., $X$ is weighted
  homogeneous. Then \cite[Satz 5]{MCrelle} applies.
  \hfill$\Box$
\end{trivlist}
From now on, $G$ will always denote the maximal reductive automorphism
group of the normal surface singularity $X$. As mentioned in the preceding
proof, $G$ is finite if $X$ is not weighted homogeneous. Otherwise,
we may assume that the
connected component $G_1$ of the unit contains at least the $\C^*$
defining the weighted homogeneous structure. In fact,
Scheja and Wiebe \cite[section 3]{SW2} and
Wahl \cite[3.6.2]{WPSPM} showed
that $G_1=\C^*$ unless $X$ is a cyclic quotient singularity.
We shall reprove this below with a different method which, at the same time,
will yield information on the finite group $G/G_1$.
\\[1.5ex]
Let $(\tilde{X},E)\to (X,0)$ be the minimal good resolution with exceptional
divisor $E$. There is a natural homomorphism $\Aut(X,0)\to\Aut(\tilde{X},E)$
obtained from the universal property of the minimal good resolution.
Here $\Aut(\tilde{X},E)$ denotes the group of germs of automorphisms
of $\tilde{X}$ along $E$.
In the weighted homogeneous case, the minimal good resolution was described
by Orlik and Wagreich \cite[Theorem 2.3.1]{OW}: If $X$ is not a cyclic
quotient singularity then $E$ has a component
$E_0$ which is uniquely determined by the property that $E_0$ has positive
genus $g$ or intersects at least three other components of $E$. The
exceptional divisor is star shaped with central curve $E_0$ and a certain
number $r$ of branches of rational curves. Here $r\ge3$ if $E_0$ is rational
and $r\ge0$ otherwise. In deriving this from \cite{OW} one has to be
aware of a result of Brieskorn \cite[Korollar 2.12]{B}: If $X$ can be
resolved by a
chain of rational curves then $X$ is a cyclic quotient singularity.
We obtain a natural homomorphism
$\Aut\: X\to\Aut\: E_0$. Recall that $\Aut\: E_0$ is finite if $E_0$ has
genus $g\ge2$ and $\Aut\: E_0=\PSL(2,\C)$ if $E_0$ is rational. Finally, if
$E_0$ is an elliptic curve then
$\Aut\: E_0=E_0\rtimes\Aut_0\, E_0$. Here the Abelian group $E_0$ acts on
itself by translations and the group $\Aut_0\, E_0$ of automorphisms
fixing the zero element is cyclic of
order 6, 4 or 2 if the $j$-invariant of $E_0$ is 0, 1 or else.
\begin{trivlist}
  \item[] \bf Theorem 2. \it
  Suppose that $X$ is weighted homogeneous but not a cyclic quotient
  singularity. Then the restriction $\rho:G\to\Aut\: E_0$
  of the natural map described above
  has kernel $\C^*$ and finite image. Hence $G_1=\C^*$ and $\rho$
  embeds $G/G_1$ into $\Aut\: E_0$.
  \item[] Proof. \rm
  Let $X\sub(\C^n,0)$ be a minimal embedding of the singularity $X$ in a
  smooth germ.
  By \cite[Satz 6]{MCrelle} the action of $G$ on $X$ can be extended to an
  action on $(\C^n,0)$ which is linear in suitable coordinates. Moreover,
  we may assume that $X$ is defined by equations which are weighted
  homogeneous polynomials in the chosen coordinates. Hence the action on
  the analytic space germ $X$ is induced from a global rational action
  on the affine algebraic variety $X$. We claim that the action on the
  germ $(\tilde{X},E)$ is induced from a global rational
  action $G\times \tilde{X}\to\tilde{X}$ on the non-singular surface
  $\tilde{X}$. In fact, one obtains a $G$-equivariant good resolution
  $X'\to X$ by succesively blowing up in $G$-invariant centres and
  normalizing.
  By the universal properties of blowing up and normalization, the
  action $G\times X'\to X'$ is rational. The minimal good resolution is
  obtained from $X'$ by blowing down $G$-invariant systems of exceptional
  curves of the first kind. This gives the
  claim. We conclude that $G\times E_0\to E_0$ is rational.
  \item[]
   The closed normal subgroup $K=\ker\:\rho\le G$ is reductive.
   Consider the representation of $K$ on
   the tangent space $\T_p\tilde{X}$ at some point $p\in E_0$. By Cartan's
   Uniqueness Theorem \cite{K} it is faithful. Since $K$ is trivial on
   the subspace $\T_p E_0$ and the representation is completely reducible
   we see that $K$ has a faithful one dimensional representation. Thus $K$
   is a subgroup of $\C^*$. It will follow that $K=\C^*$ if we
   show that $\im\:\rho$ is finite.
   This is obvious if $E_0$ has genus $g\ge2$. But also for $g=1$ since
   $\C^*$ and hence every connected reductive group can only act trivially
   on an elliptic curve. In the remaining
   case $g=0$ consider the $r$ intersection
   points of $E_0$ with other components of $E$. They have to be
   permuted by any element of $\im\:\rho$. We obtain a homomorphism
   $\im\:\rho\to{\rm S}_r$ to the symmetric group which is injective since
   $r\ge3$ and since automorphisms of the projective line have at most
   two fixed points.
   \hfill $\Box$
\end{trivlist}
\begin{trivlist}
  \item[] {\it Remarks.} \rm
  (i) Let $\G$ denote the weighted dual graph defined by the minimal
  good resolution and $\Aut\:\G$ its group of automorphisms. Let $r$ be
  the number of branches emanating from the central curve of genus $g$.
  We have a natural homomorphism $\im\:\rho\to\Aut\:\G$. A non-trivial
  automorphism of $E_0$ has at most $2g+2$ fixed points, \cite[V.1.1]{FK}.
  Thus, if
  $r>2g+2$ then $G/G_1$ embeds into $\Aut\:\G$. This is always the case
  if $g=0$ since then $r\ge3$. The latter result was previously obtained
  (with a different proof) by Aumann-K\"orber \cite[3.9]{A-K}.
  \item[] (ii) If $G/G_1\to\Aut\:\G$ is injective we obtain sufficient
  conditions in terms of the weighted resolution graph for $G/G_1$ to be
  trivial or cyclic. (Recall \cite[III.7.7]{FK}
  that a finite group of automorphisms of a smooth
  complete curve fixing a point is cyclic.)
  \item[] (iii) The exact sequence
  $1\to G_1\to G\to G/G_1\to 1$ splits if and only if
  $G\simeq G_1\times G/G_1$ is a direct product over $G_1$. This follows
  from the fact \cite[Proposition 3.10]{WPSPM} that $G_1=\C^*$ is central
  in $G$.
  \item[] (iv) Let $X=\C^2/H$ be a cyclic quotient singularity where the
  group $H$ is generated by
  $\left(\begin{array}{cc}
                       \zeta & 0 \\
                       0     & \zeta^q   \end{array}\right)$
  with $\zeta^n=1$ and $(q,n)=1$. Then $G_1=\GL(2,\C)/H$ if $q=1$ and
  $G_1=\C^*\times\C^*$ else, see \cite[3.6.2]{WPSPM}.
  Also in this case, there is a natural homomorphism $G/G_1\to\Aut\:\G$
  which, in fact, is bijective, \cite[3.11, 3.12]{A-K}.
  \item[] (v) Suppose that $X$ is not necessarily weighted homogeneous
  but the exceptional
  divisor $E$ has a component $E_1$ which is fixed by every automorphism
  of $(\tilde{X},E)$. (This is, e.~g., the case for determinantal
  rational singularities of multiplicity $m\ge3$ since their exceptional
  divisor has a unique component $E_1$ of self intersection number $-m$,
  \cite[Theorem 3.4]{Wrat}.)
  Then there is again a natural homomorphism
  $G\to\Aut\: E_1$. The proof of Theorem 2 shows that its kernel is
  cyclic if $X$ is assumed to be not weighted homogeneous.
  \item[] (vi) If $X$ is rational or minimally elliptic
  but not weighted homogeneous and not a cusp singularity then the
  kernel $K$ of the natural homomorphism $G\to\Aut\:\G$ is cyclic. In
  fact, by \cite[Lemma 1.3]{B} and \cite[Proposition 3.5]{Lauell}
  all components of $E$ are rational but $\G$ is not a chain nor a cycle.
  (Note that, besides the cusp singularities, also simple elliptic and
  cyclic quotient singularities are excluded because they are weighted
  homogeneous.) Hence there is a component
  $E_1$ intersecting at least three other components. As $E_1$ is rational
  $K$ acts trivially on $E_1$. The proof of Theorem 2 shows that $K$ is
  cyclic. The special case of non weighted homogeneous exceptional
  unimodal singularities will be discussed in more detail in Example 3
  of section 4.
  For cusp singularities the kernel of $G\to\Aut\:\G$ is Abelian
  but not necessarily cyclic, \cite{Wcus}.
  \item[] (vii) When we are considering the maximal reductive automorphism
  group $G$ of some weighted homogeneous singularity $X$ we may choose
  coordinates such that, at the same time, $G$ acts linearly on the ambient
  space and $X$ is defined by weighted homogeneous polynomials, see the
  beginning of the proof of Theorem 2. This does not mean that if we start
  with a weighted homogeneous polynomial defining some $X$ then $G$ will
  be linear in the given coordinates. For instance, let $X\sub(\C^3,0)$
  be defined by $f=x_1^3-2x_3x_2^2-x_3^2$ which is weighted homogeneous of
  weights 4, 3, 6 and degree 12. The $\C^*$-action is given by
  $t\cdot(x_1,x_2,x_3)=(t^4x_1,t^3x_2,t^6x_3)$. Suppose that
  $G\le\Aut(\C^3,0)$
  contains this $\C^*$ and defines a maximal reductive subgroup of
  $\Aut\: X$. Since $\C^*$ is central in $G$, see Remark (iii) above,
  it is easily seen that $G\cap\GL(3,\C)=\C^*$. But
  $G$ is larger than $\C^*$ since $X$ is isomorphic to the singularity
  $E_6$ defined by $g=x_1^3+x_2^4-x_3^2$ whose maximal reductive automorphism
  group clearly contains $\C^*\times\Z_2$ where $\Z_2$ acts by
  $(x_1,x_2,x_3)\mapsto(x_1,x_2,-x_3)$. It follows from Remark (i) that,
  in fact, $\C^*\times\Z_2$ is the maximal reductive automorphism group of
  $E_6$ and hence of $X$.
\end{trivlist}
\section{The finite group $G/G_1$}
In the weighted homogeneous case,
we are going to study the finite group $G/G_1$ viewed as a subgroup of
the automorphism group of the central curve via the natural homomorphism
$\rho$.
First consider singularities $X$ such that the exceptional divisor is
irreducible, $E=E_0$, of genus $g$. Let $N\to E_0$ be the normal bundle
of $E_0\sub\tilde{X}$ with zero-section $E_0\sub N$. It follows from work
of Grauert \cite[\S 4]{Grau} (see \cite[6.2]{Wag}) that the germs
$(\tilde{X},E_0)$ and $(N,E_0)$ are isomorphic if the self intersection number
satisfies $E_0\cdot E_0<4 - 4g$. We determine $G$ for these singularities:
\begin{trivlist}
  \item[] \bf Theorem 3. \it
  Let $\pi:N\to E_0$ be a negative line bundle on a smooth complete curve
  $E_0$ of genus $g\ge1$, and let $X$ be the singularity obtained by
  contracting the zero-section $E_0\sub N$.
  \item[] {\rm (i)} Let $\tilde{G}$ be the group of automorphisms $g$ of the
  manifold $N$ which restrict to an automorphism $g_0$ of the zero-section
  such that $g_0\circ\pi=\pi\circ g$. Then $\tilde{G}$ is the maximal
  reductive automorphism group of $X$.
  \item[] {\rm (ii)} Let $g=1$ and $-b=E_0\cdot E_0$. Then the natural map
  $\rho:G\to\Aut\: E_0$
  has image $(\Z_b\times\Z_b)\rtimes\Aut_0\, E_0$ where
  $\Z_b\times\Z_b$ is the group of $b$-torsion points of the Abelian group
  $E_0$ and $\Aut_0\, E_0$ is the group of automorphisms fixing the zero
  element. The exact sequence $1\to G_1\to G\to G/G_1\to 1$ splits
  if and only if $b=1$.
  \item[] Proof. \rm
  Let $\tilde{\rho}:\tilde{G}\to\Aut\: E_0$ be the restriction map. It has
  kernel $\C^*$. Let $g=1$. We may assume that $N=\O(-b\cdot p_0)$ for
  some $p_0\in E_0$. Further we may assume that $p_0$ is the zero element $0$
  of the group $E_0$. We claim that $\phi\in\Aut\: E_0$ belongs to
  $\im\:\tilde{\rho}$ if and only if $p=\phi(0)$ is a $b$-torsion point
  of $E_0$. In fact, for $\psi=\phi^{-1}$ let
  $\psi_*:\O(-b\cdot p)\to\O(-b\cdot 0)$ be the induced map satisfying
  $\psi\circ\pi=\pi\circ\psi_*$ (where $\pi$ denotes both bundle
  projections). If $\phi=g_0=\tilde{\rho}(g)$ for some $g\in\tilde{G}$
  then $g\circ\psi_*:\O(-b\cdot p)\to\O(-b\cdot 0)$ is an isomorphism of
  line bundles. Conversely, if $\Psi:\O(-b\cdot p)\to\O(-b\cdot 0)$ is a
  line bundle isomorphism then $\Psi\circ \phi_*$ is an element of $\tilde{G}$
  restricting to $\phi$. Now observe that the line bundles $\O(-b\cdot p)$
  and $\O(-b\cdot 0)$ are isomorphic if and only if the divisors
  $-b\cdot p$ and $-b\cdot 0$ are linearly equivalent if and only if
  $bp=0$ in the group $E_0$. This proves the claim. It follows from
  $\Aut\: E_0= E_0\rtimes\Aut_0\,E_0$ that
  $\im\:\tilde{\rho}=(\Z_b\times\Z_b)\rtimes\Aut_0\, E_0$. Since
  $\tilde{\rho}$ has finite image also for $g\ge2$ we conclude that
  $\tilde{G}$ is reductive in any case.
  \item[]
  Now $\tilde{G}$ induces a reductive algebraic subgroup of $\Aut\: X$.
  Hence $\tilde{G}$ is contained in a maximal reductive subgroup $G$
  viewed as a group of automorphisms of the germ $(N,E_0)$. As both
  $\tilde{G}$ and $G$ are one dimensional they are equal if they have the
  same image in $\Aut\: E_0$. But every $g\in G$ induces a
  $\tilde{g}\in\tilde{G}$ with the same restriction to $E_0$ since
  $N\to E_0$ is the normal bundle of $E_0\sub N$.
  \item[]
  Let us return to the case $g=1$. If $b=1$ then the divisor $-0$ defining
  $N$ is invariant under $\im\:\rho=\Aut_0\, E_0$. This yields a section
  $\im\:\rho\to\tilde{G}=G$.
  To prove the non-existence of such a section for $b\ge2$
  we write $E_0=\C/\L$ where the lattice
  $\L$ is generated by the primitive periods $\w_1,\w_2$ with
  ${\rm Im}\:\w_2/\w_1>0$. Let $\phi:E_0\to E_0$ be the translation
  given by $z\mapsto z+\a$ with $\a=\w_1/b$. Set $\a_0=-(b+1)/2b\cdot\w_1$
  and $\a_k=\a_0+k\a$ for
  $k\in\Z$. Then the divisor $D_1=\sum_{k=1}^b\a_k$ on $E_0$ is
  $\phi$-invariant and $\phi$ induces an automorphism of $\O(-D_1)$ over
  $\phi$. It is given by
  \[s_1(z)\mapsto s_1(z+\a)
  \]
  with $s_1$ denoting a rational section of $\O(-D_1)$. Because
  $\a_1+\ldots+\a_b=0$ in $\C$ the divisors $D_1$ and
  $D_0=b\cdot 0$ are linearly equivalent, say $D_1-D_0=(h_1)$ for some
  rational function $h_1$ on $E_0$. Using a rational section $s_0$ of
  $N=\O(-D_0)$ one defines a line bundle isomorphism $N\to\O(-D_1)$ by
  \[s_0(z)\mapsto h_1(z)s_1(z).
  \]
  We obtain an automorphism $\tilde{\phi}\in\tilde{G}=G$ of $N$ over $\phi$
  given by
  \[s_0(z)\mapsto\frac{h_1(z)}{h_1(z+\a)}s_0(z+\a).
  \]
  Now let $\b=\w_2/b$, $\psi$ the corresponding translation,
  $\b_0=-(b+1)/2b\cdot\w_2$, $\b_k=\b_0+k\b$,
  $D_2=\sum_{k=1}^b\b_k$, $D_2-D_0=(h_2)$ and $\tilde{\psi}\in G$
  the induced map. Then
  $\tilde{\psi}^{-1}\tilde{\phi}^{-1}\tilde{\psi}\tilde{\phi}\in G$
  is given by
  \[s_0(z)\mapsto\frac{h(z)h(z+\a+\b)}{h(z+\a)h(z+\b)}s_0(z)
  \]
  where $h=h_1/h_2$. As this commutator is contained in
  the kernel of $\rho$ the function
  $h(z)h(z+\a+\b)h(z+\a)^{-1}h(z+\b)^{-1}$ is constant, say equal to
  $\l\in\C^*$. If we assume that $\rho$ admits a section, hence
  that $G=\C^*\times B$ with a subgroup $B$, see Remark (iii) of section 2,
  then the commutator must be contained in $B$ and $\l=1$. We
  are going to show that $\l=\exp(2\pi i/b)$, hence $b=1$.
  \item[]
  Let $\zeta$ and $\s$ be Weierstrass' $\zeta$- and $\s$-function.
  One has
  \[\s(z+\w_i)=-\exp(\eta_i(z+\w_i/2))\cdot\s(z)
  \]
  with $\eta_i=2\zeta(\w_i/2)$ for $i=1,2$, see
  \cite[II.1.13, Satz 3]{HC}. The elliptic function $h$ has
  $\a_1,\ldots,\a_b$ and $\b_1,\ldots,\b_b$ as complete systems of zeros
  and poles. Since $\a_1+\ldots+\a_b=\b_1+\ldots+\b_b$ it follows from
  \cite[II.1.14, Satz 1]{HC} that
  \[h(z)=\prod_{k=1}^b\frac{\s(z-\a_k)}{\s(z-\b_k)}
  \]
  up to a constant. One calculates
  \begin{eqnarray*}
  \frac{h(z)h(z+\a+\b)}{h(z+\a)h(z+\b)} & = &
  \frac{\s(z-\a_b)\s(z-\b_0)\s(z+\b-\a_0)\s(z+\a-\b_b)}
       {\s(z-\a_0)\s(z-\b_b)\s(z+\b-\a_b)\s(z+\a-\b_0)} \\
  & = & \exp(\eta_1\b-\eta_2\a) \\
  & = & \exp((\eta_1\w_2-\eta_2\w_1)/b).
  \end{eqnarray*}
  Legendre's relation $\eta_1\w_2-\eta_2\w_1=2\pi i$, see
  \cite[II.1.11]{HC}, gives the claim.
  \hfill$\Box$
\end{trivlist}
Let $X$ be weighted homogeneous but not a cyclic quotient singularity.
For $i=1,\ldots,r$ let $E_{ij}$, $j=1,\ldots,l_i$, be the curves on
the $i$-th branch
of the exceptional divisor of the minimal good resolution. Assume that
$E_{i1}$ intersects the central curve $E_0$, say in the point $p_i$,
and that $E_{ij}$ intersects $E_{i,j+1}$. Let $-b=E_0\cdot E_0$ and
$-b_{ij}=E_{ij}\cdot E_{ij}$ denote the self intersection numbers. Finally,
consider the Hirzebruch-Jung continued fractions
\[\a_i/\b_i=b_{i1}-1/(b_{i2}-1/(\ldots-1/b_{i,l_i})\ldots)
\]
with coprime positive integers $\b_i<\a_i$.
\begin{trivlist}
  \item[] \bf Theorem 4. \it
  Let $X$ be weighted homogeneous but not a cyclic quotient singularity
  and suppose that the central curve is rational.
  \item[] {\rm (i)} Then the image of $\rho:G\to\Aut\: E_0$ equals
  the group $A$ of automorphisms of $E_0$ which permute the points
  $p_1,\ldots,p_r$ while preserving the pairs $(\a_i,\b_i)$.
  \item[] {\rm (ii)} If $b$ is even or if $A$ is cyclic or if $A$ is
  dihedral of order $2q$ with $q$ odd then the maximal reductive automorphism
  group is a direct product $G\simeq\C^*\times A$.
  \item[] Proof. \rm
  Obviously, $\im\:\rho\sub A$. For the other inclusion
  we use Pinkham's description \cite[Theorem 5.1]{P} of the affine
  coordinate ring of a weighted homogeneous surface. Let the divisor $D_0$
  correspond to the conormal bundle of $E_0\sub\tilde{X}$. Consider the
  divisor
  \[D=D_0-\sum_{i=1}^r\b_i/\a_i\cdot p_i
  \]
  with rational coefficients.
  Then the coordinate ring of the affine algebraic variety $X$ is isomorphic,
  as a graded ring,
  to $\bigoplus_{k=0}^\infty L(kD)$ where $L(kD)$ denotes the space
  of rational functions $f$ on $E_0$ with $(f)\ge-kD$. Now recall the
  classification of the finite subgroups of $\PSL(2,\C)$ and their orbits
  on the projective line, \cite[Chapter I, \S 6]{L}. The cyclic groups
  have a fixed point. The dihedral group of order $2q$ has orbits of
  length 2 and $q$. And for the tetrahedral, the octahedral and the
  icosahedral group the greatest common divisor of the lengths of the
  orbits is 2. Thus, if $E_0$ is rational and one of the hypotheses of (ii)
  is fulfilled then there exists on $E_0$ an $A$-invariant divisor $D_1$
  of degree $b$. The sum in the definition of $D$ is $A$-invariant. Since
  $E_0$ is rational the two divisors $D_0,D_1$ of the same degree $b$
  are linearly equivalent. Hence we may replace $D_0$ by $D_1$ and assume
  that $D$ itself is $A$-invariant.
  Now there is an obvious action of $A$ on
  $\bigoplus_{k=0}^\infty L(kD)$. We obtain a homomorphism
  $A\to\Aut\: X:\phi\mapsto\phi'$ whose image may be assumed to lie in $G$.
  One checks that $\rho(\phi')=\phi$ for all $\phi\in A$. Hence
  $A=\im\:\rho$ and the exact sequence
  $1\to\C^*\to G\to A\to 1$ splits. Then $G\simeq \C^*\times A$, see
  Remark (iii) of section 2.
  To prove (i) in the cases not covered by
  (ii) apply (ii) to the cyclic groups generated by the elements of $A$.
  \hfill$\Box$
\end{trivlist}
\begin{trivlist}
  \item[] \bf Corollary. \it
  Let $X$ be weighted homogeneous such that the exceptional divisor
  consists of exactly three branches emanating from a rational central curve.
  Then $G\simeq\C^* \times \Aut\:\G$ where $\G$ denotes the weighted dual
  resolution graph.
  \item[] Proof. \rm
  Recall from Remark (i) of section 2 the injection $A\to\Aut\:\G$. In the
  present case it is an isomorphism since three points on the projective line
  have no moduli. Part (ii) of the Theorem applies as $\Aut\:\G$ is
  trivial or cyclic of order 2 or isomorphic to the symmetric group
  ${\rm S}_3$, i.\ e., dihedral of order 6.
  \hfill$\Box$
\end{trivlist}
\begin{trivlist}
  \item[] \it Remarks. \rm
  (i) The Corollary applies to all non-cyclic quotient singularities.
  In this special case the result was previously obtained by
  Aumann-K\"orber \cite[3.12]{A-K}. She uses the fact \cite[3.6.3]{WPSPM}
  that for a quotient singularity $\C^2/H$ the maximal reductive
  automorphism group is ${\rm N}(H)/H$ where ${\rm N}(H)$ denotes the
  normalizer of $H$ in $\GL(2,\C)$. Then she explicitly computes this
  normalizer for every $H$. The Corollary applies, as well, to triangle
  singularities, see Example 2 below.
  \item[] (ii) The proof of Theorem 4 shows $G\simeq\C^*\times A$ if $g$ is
  arbitrary and the divisor $\sum_{i=1}^r p_i$ corresponds to the
  conormal bundle. Such singularities (with $g=1$) have been considered
  by Tomaru \cite{T}.
\end{trivlist}
We end this section by showing that every finite group appears
as $G/G_1$ for a suitable $X$.
\begin{trivlist}
  \item[] \bf Theorem 5. \it
  Let $E_0$ be an arbitrary smooth complete curve and let $A\le\Aut\: E_0$
  be an arbitrary finite subgroup. Then there is a weighted homogeneous
  normal surface singularity $X$ with central curve $E_0$ and maximal
  reductive automorphism group $G\simeq\C^*\times A$.
  \item[] Remark. \rm Hurwitz showed that every finite group
  can be realized as a subgroup of $\Aut\: E_0$ for some smooth complete
  curve $E_0$, see \cite{Green}.
  \item[] {\it Proof of the Theorem.}
  We first show that there exist finitely many $A$-orbits $B_1,\ldots,B_k
  \sub E_0$ such that $A$ consists of exactly those
  $\phi\in\Aut\: E_0$ with $\phi(B_j)=B_j$ for all $j=1,\ldots,k$.
  Take some $A$-orbits $B_1,\ldots,B_l$ and let $A_1$ be the group of
  $\phi\in\Aut\: E_0$ with $\phi(B_j)=B_j$ for $j=1,\ldots,l$. If the
  number of elements of $B_1\cup\ldots\cup B_l$ is sufficiently big
  (say, $\ge3$ if $E_0$ has genus $g=0$ and $\ge5$ if $g=1$) then $A_1$
  is finite. We have $A\sub A_1$. If the inclusion is strict then a generic
  $A_1$-orbit (not containing any of the finitely many points which are
  fixed by some element of $A_1$) has more elements than $A$. Choose an
  $A$-orbit $B_{l+1}$ contained in a generic $A_1$-orbit. Then the group
  $A_2$ of elements of $A_1$ which moreover map $B_{l+1}$ onto itself
  satisfies $A\sub A_2\subsetneq A_1$. After finitely many steps we arrive
  at the claim. Now attach pairs $(\a_i,\b_i)$ of coprime positive integers
  $\b_i<\a_i$ to the points $p_1,\ldots,p_r$ of
  $B_1\cup\ldots\cup B_k$ in such a way that two points get the same label
  if and only if they belong to the same orbit.
  Then $A$ is the group of
  automorphisms of $E_0$ which permute the points $p_1,\ldots,p_r$
  while preserving the pairs $(\a_i,\b_i)$. Moreover, choose an $A$-invariant
  divisor $D_0$ on $E_0$ of degree $b>\sum_{i=1}^r \b_i/\a_i$. Then
  Pinkham's construction yields a weighted homogeneous $X$ with the given
  data and such that its maximal reductive automorphism group is isomorphic
  to $\C^*\times A$, see the proof of Theorem 4.
  \hfill$\Box$
\end{trivlist}
\section{Examples}
Catanese \cite[section 2]{C} has studied involutions of rational double
points.
We extend this to all weighted homogeneous singularities with rational
central curve and three branches.
To determine the occuring quotients we need two Lemmas. As usual, the weight
$-2$ is omitted in the resolution graphs.
\begin{trivlist}
  \item[] \bf Lemma 1. \it
  Let $X=X_{b,1}$ be the cyclic quotient singularity with weighted graph
  \[\Aone{-b}
  \]
  Let $\s$ be an involution of $X$ and consider its action on the minimal
  resolution $\tilde{X}$ and on the exceptional curve $E$.
  \item[] {\rm (i)} If $\s$ fixes $E$ pointwise then
  $X/\s$ has weighted graph
  \[\Aone{-2b}
  \]
  Otherwise there are the following cases:
  \item[] {\rm (ii)} If $b$ is odd then one of the two fixed points of $\s$
  on $E$ is an isolated fixed point on $\tilde{X}$ and the other is not.
  Then $X/\s$ has weighted graph
  \[\Atwo{-(b+1)/2}{}
  \]
  \item[] {\rm (iii)} If $b$ is even then either both fixed points of $\s$
  on $E$ are isolated fixed points on $\tilde{X}$ or both are not. The
  corresponding weighted graphs are
  \[\Athree{}{-(b+2)/2}{}
  \qquad\mbox{and}\qquad
  \Aone{-b/2}
  \]
  \item[] Proof. \rm
  Let $q$ be a fixed point of $\s$ lying on $E$. If $q$ is an isolated
  fixed point it gives rise to an $A_1$-singularity in $\tilde{X}/\s$
  which can be resolved by inserting a rational $(-2)$-curve. Otherwise,
  $q$ is mapped onto a smooth point of $\tilde{X}/\s$. A resolution of
  $\tilde{X}/\s$ clearly resolves $X/\s$. Thus, according to
  whether there are two or one or no isolated fixed points lying on $E$
  the quotient $X/\s$ has weighted graph
  \[\Athree{}{-c}{}
  \qquad\mbox{or}\qquad
  \Atwo{-c}{}
  \qquad\mbox{or}\qquad
  \Aone{-c}
  \]
  with some $c$. To determine $c$ and to prove the other assertions
  write $X=X_{b,1}=\C^2/H$ where $H$ is generated by
  $\left(\begin{array}{cc}
                       \zeta & 0 \\
                       0     & \zeta   \end{array}\right)$
  and $\zeta$ is a primitive $b$-th root of unity. As mentioned in Remark
  (iii) of section 2 the maximal reductive automorphism group of $X$ is
  $\GL(2,\C)/H$. Hence $X/\s\simeq\C^2/\Sigma$ with some group
  $\Sigma\le\GL(2,\C)$ containing $H$ as a subgroup of index two.
  Clearly $\Sigma $ is Abelian. So we may assume that it consists of
  diagonal matrices.
  \item[]
  Consider first the case that $\Sigma$ is cyclic, say generated by $\t$.
  We may assume that $\t^2=
  \left(\begin{array}{cc}
                       \zeta & 0 \\
                       0     & \zeta   \end{array}\right)$.
  If $\t=
  \left(\begin{array}{cc}
                       \eta & 0 \\
                       0     & \eta   \end{array}\right)$
  with a $2b$-th root of unity $\eta$ then $X/\s\simeq\C^2/\t=X_{2b,1}$ with
  graph as in (i). On
  $\tilde{X}=\{(z_1,z_2,(w_1:w_2))\in \C^2\times\P_1,\; z_1w_2^b=z_2w_1^b\}$
  we have $\s(z,w)=(-z,w)$ and $\s$ fixes $E=0\times\P_1$ pointwise. The
  next possibility is $\t=
  \left(\begin{array}{cc}
                       \eta & 0 \\
                       0     & -\eta   \end{array}\right) =
  \left(\begin{array}{cc}
                       \eta & 0 \\
                       0     & \eta^{b+1}   \end{array}\right)$
  where again $\eta$ is a $2b$-th root of unity. If $b$ is even then $2b$
  and $b+1$ are coprime and $X/\s\simeq\C^2/\t=X_{2b,b+1}$. The continued
  fraction expansion of $2b/(b+1)$ shows that the graph is the first one
  in (iii). If $b$ is odd then $b$ and $(b+1)/2$ are coprime and
  $\t^b=
  \left(\begin{array}{cc}
                       -1 & 0 \\
                       0     & 1   \end{array}\right)$.
  Since $\t$ acts on $\C^2/\t^b\simeq\C^2$ by
  $\left(\begin{array}{cc}
                       \eta^2 & 0 \\
                       0     & \eta^{b+1}   \end{array}\right) =
  \left(\begin{array}{cc}
                       \zeta & 0 \\
                       0     & \zeta^{(b+1)/2}   \end{array}\right)$
  we see $X/\s\simeq X_{b,(b+1)/2}$ with graph as in (ii). In these two cases
  the assertion on the fixed points is obvious from the graphs.
  \item[]
  Now consider the case that $\Sigma$ is not cyclic. Then
  $\Sigma=H\times<\t>$ with some involution $\t\in\GL(2,\C)$, and $b$ must
  be even. As
  $\left(\begin{array}{cc}
                       -1 & 0 \\
                       0     & -1   \end{array}\right)\in H$ we have
  $\t=\left(\begin{array}{cc}
                       -1 & 0 \\
                       0     & 1   \end{array}\right)$ without loss of
  generality and $\Sigma$ contains the reflection group $T$ generated by
  $\t$ and $-\t$.
  The generator of $H$ acts on $\C^2/T\simeq\C^2$ by
  $\left(\begin{array}{cc}
                       \zeta^2 & 0 \\
                       0     & \zeta^2   \end{array}\right)$.
  Hence $X/\s\simeq X_{b/2,1}$ and the graph is the second one in (iii).
  Moreover, $\s(z,w)=(z,-w)$ on $\tilde{X}$.
  \hfill$\Box$
\end{trivlist}
\begin{trivlist}
  \item[] {\it Remark.} The statement of the Lemma is true also for $b=1$,
  i.~e., if
  $X=(\C^2,0)$ and $\tilde{X}$ is the blow up of $0$. For $b=2$ it
  appears in \cite[p.~80]{S}.
\end{trivlist}
\begin{trivlist}
  \item[] \bf Lemma 2. \it
  Let $X$ be weighted homogeneous but not a cyclic quotient singularity.
  Let $\phi$ be an automorphism of $X$ of finite order and consider its
  action on the minimal good resolution $\tilde{X}$. Suppose that $\phi$
  fixes the intersection point $p$ of two components $E_1,E_2$ of the
  exceptional divisor $E$ and that $p$ is not an isolated fixed point of
  $\phi$ on $\tilde{X}$. Then near $p$ the fixed point locus coincides with
  $E_1$ or $E_2$.
  \item[] Proof. \rm
  Near $p$ the automorphism $\phi$ can be linearized. In suitable local
  coordinates $x,y$ it is given by $\phi(x,y)=(x,\l y)$ with some $\l\not=1$.
  As $E$ is star shaped it is not possible that $\phi$ interchanges $E_1$
  and $E_2$. Hence they are left invariant. It is easily seen that a
  smooth $\phi$-invariant curve through $p$ different from the fixed point
  locus $\{y=0\}$ must be tangent to $\{x=0\}$. Since $E_1$ and $E_2$
  intersect transversely the Lemma is proven.
  \hfill$\Box$
\end{trivlist}
\begin{trivlist}
  \item[] \bf Proposition 1. \it
  Let $X$ be weighted homogeneous such that the exceptional divisor consists
  of exactly three branches emanating from a rational central curve.
  \item[] {\rm (ii)} If $\Aut\:\G$ is trivial then $X$ has (up to
  conjugation) exactly one involution, namely the one contained in $\C^*$.
  \item[] {\rm (ii)} If $\Aut\:\G$ is not trivial then besides
  $\s_1\in\C^*$ there are (up to conjugation) exactly two more involutions
  $\s_2$ and $\s_3=\s_1\s_2$ in $\Aut\: X$. They can be distinguished by
  the property that for one of them, say $\s_2$, the quotient $X/\s_2$
  is a cyclic quotient singularity whereas $X/\s_3$ is not. Here we agree
  that the smooth germ $(\C^2,0)$ is called a cyclic quotient singularity,
  too.
  \item[] Proof. \rm
  (i) is clear from the Corollary of Theorem 4. In (ii) we have
  $G\simeq\C^*\times \Aut\:\G$ with $\Aut\:\G=\Z_2$ or $\rm S_3$. Since
  $\rm S_3$
  has, up to conjugation, only one involution there are only three
  involutions to consider in $G$. We may assume that $\s_2$ and $\s_3$
  interchange the second and third branch which therefore get identified
  in the quotient. Both involutions fix the intersection point of $E_0$
  with the first branch and have a second fixed point $q$ on the
  projective line $E_0$. In suitable local coordinates $x,y$ around $q$
  with $E_0=\{x=0\}$ we have $\s_1(x,y)=(-x,y)$, $\s_2(x,y)=(x,-y)$ and
  $\s_3(x,y)=(-x,-y)$. Hence $q$
  is mapped onto a smooth point in $\tilde{X}/\s_2$, but in the resolution
  of $\tilde{X}/\s_3$ there appears a new branch consisting of a single
  $(-2)$-curve.
  The only isolated fixed points possibly lying on the first branch are
  intersection points of components plus, maybe, one point on the curve
  at the end of the branch. Thus, after resolving the $A_1$-singularities
  appearing in the quotient we are left with a chain of rational curves.
  We conclude that $X/\s_2$ is a cyclic quotient singularity. And it will
  follow that $X/\s_3$ has a star shaped graph with three branches (and
  hence is not a cyclic quotient singularity) if we can show that the
  chain of rational curves arising from the first branch of $X$ cannot be
  blown down completely to a smooth point. From Lemma 1 one sees that there
  is only one possible way for a $(-1)$-curve to appear. Namely, if the
  branch contains a $(-2)$-curve $E_i$, not fixed pointwise, such that the
  two fixed points on $E_i$ are not isolated fixed points of the involution
  on $\tilde{X}$. But then Lemma 2 shows that the neighbouring components
  have to be fixed pointwise. Using Lemma 1 again we see that the
  configuration
  \[\Athree{-b_{i-1}}{}{-b_{i+1}}
  \]
  produces
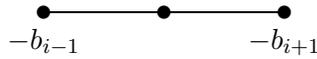
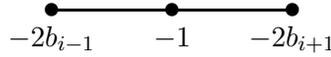
  \[\Athree{-2b_{i-1}}{-1}{-2b_{i+1}}
  \]
  in the quotient, hence
  \[
  \unitlength0.8mm
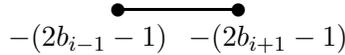
  \begin{picture}(30,9)
  \multiput(5,6)(20,0){2}{\circle*{2}}
  \put(0,0){\makebox(0,0)[b]{$-(2b_{i-1}-1)$}}
  \put(30,0){\makebox(0,0)[b]{$-(2b_{i+1}-1)$}}
  \put(5,6){\line(1,0){20}}
  \end{picture}
  \]
  after blowing down the $(-1)$-curve. Consequently, the blowing down does
  not create new $(-1)$-curves. This implies the claim.
  \hfill$\Box$
\end{trivlist}
\begin{trivlist}
  \item[] {\it Example 1.} Proposition 1 applies to each non-cyclic
  quotient singularity $X$.
  As can be seen from \cite{R} the exceptional divisor has a branch
  consisting of a single $(-2)$-curve. It follows from Lemma 1 that this
  branch disappears in the quotient with respect to the involution
  $\s_1\in\C^*$. As in the proof of Proposition 1 one then shows that
  $X/\s_1$ is a cyclic quotient singularity. Of course, it is a simple
  task to determine it explicitly in each case. From \cite{R} one sees
  that $\Aut\:\G$ is non-trivial if and only if $X$ is dihedral or
  $X=T_m$ is tetrahedral with $m\equiv1$ or 5 mod 6. In these cases let
  $\s_3$ be the involution for which $X/\s_3$ is not a cyclic quotient
  singularity. If $X$ is dihedral then $X/\s_3$ is dihedral again. And
  for $X=T_m$, $m\equiv1$ or 5 mod 6, one obtains $T_m/\s_3\simeq O_m$, an
  octahedral singularity. More precisely, if $m=6(b-2)+1$ then $T_m$ has
  graph
  \[\Esix{}{}{}{-b}{}{}
  \]
  Using Lemmas 1 and 2 one sees that $T_m/\s_3$ has graph
  \[\Dfive{}{-4}{-(b+1)/2}{}{}
  \]
  if $b$ is odd, but
  \[\Eseven{}{}{}{-(b+2)/2}{}{}{}
  \]
  if $b$ is even. For $m=6(b-2)+5$ the graphs of $T_m$ and
  $O_m\simeq T_m/\s_3$
  are
  \[\Dfour{}{-3}{-b}{-3}
  \]
  and
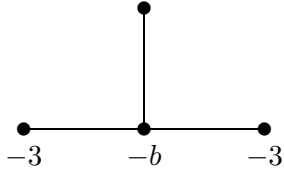
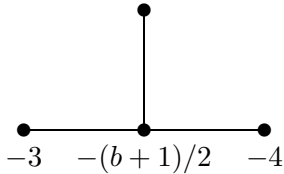
  \[\Dfour{}{-3}{-(b+1)/2}{-4}
  \]
  if $b$ is odd, but
  \[\Dsix{}{-3}{-(b+2)/2}{}{}{}
  \]
  if $b$ is even.
\end{trivlist}
Before discussing the next examples we make a digression.
\begin{trivlist}
  \item[] \bf Lemma 3. \it For a normal surface singularity $X$ let
  \[V(X)=\G(X-0,\Omega^2)/L^2(X-0)
  \]
  be the vector space of $2$-forms defined
  on a deleted neighbourhood of the singular point modulo the subspace
  of square integrable forms. Let $H\le\Aut\:X$ be a finite subgroup and
  $\pi:X\to X/H$ the quotient map. Then the pullback of forms induces an
  injection $V(X/H)\hookrightarrow V(X)^H$.
  Here the upper index $H$ denotes the subspace of $H$-invariants. For
  the geometric genus $p_g(X)=\dim\:V(X)$, see \cite[Theorem 3.4]{Laurat},
  this yields
  \[p_g(X/H)\le p_g(X)
  \]
  with strict inequality if the representation of $H$ on $V(X)$ is not
  trivial. If $\pi$ is unramified outside the singular point then
  $V(X/H)=V(X)^H$.
  \item[] Proof. \rm
  For $\a\in\G(X/H-0,\Omega^2)$ the pullback $\pi^*\a\in\G(X-0,\Omega^2)$
  is $H$-invariant. And it is square integrable if and only if $\a$ is so.
  If $\pi$ is unramified outside $0$ then every $H$-invariant form $\a'$ on
  $X-0$ induces a form $\a$ on $X/H-0$ with $\a'=\pi^*\a$.
  \hspace*{\fill}$\Box$
\end{trivlist}
Let $X$ be weighted
homogeneous but not a cyclic quotient singularity. Suppose that $X$ is
Gorenstein with nowhere zero $2$-form $\omega\in\G(X-0,\Omega^2)$. We
may assume \cite[p.~56]{Ga} that $\omega$ is $G$-equivariant:
$g^*\omega=\chi(g)\cdot\omega$ for all $g\in G$ with some character
$\chi:G\to\C^*$.
Let $\bar{G}$ be the kernel of $\chi$.
Ganter \cite[Lemma 8.3]{Ga} has shown
that $X\to X/\bar{G}$ is unramified outside the singular point.
In particular, $\omega$ induces on $X/\bar{G}-0$ a nowhere zero $2$-form
and $X/\bar{G}$ is Gorenstein.
There is an integer $\ep$ such that $\chi(t)=t^{-\ep}$
for all $t\in\C^*=G_1$. It is known \cite[Corollary 3.3]{Wjac} that $\ep<0$
if and only
if $\ep=-1$ if and only if $X$ is a rational double point. And $\ep=0$
if and only if $X$ is simple elliptic.
It may be mentioned at this place that Ganter \cite[Theorem 8.5]{Ga}
has obtained the bound
\[|G/G_1|\le 42\cdot(-P_X\cdot P_X)/\ep
\]
if $X$ is not log-canonical. The invariant $-P_X\cdot P_X$ can be
calculated (using the notation of section 3) as
\[-P_X\cdot P_X=
\frac{(2g-2+r-\sum_{i=1}^r1/\a_i)^2}{b-\sum_{i=1}^r\b_i/\a_i},
\]
see \cite[Theorem 3.2]{WPP}.
If we now assume that $\ep=\pm1$ then
$G=\C^*\times\bar{G}$. Thus we obtain a ``natural'' section
$G/G_1\hookrightarrow G$ of $G\to G/G_1$.
\begin{trivlist}
  \item[] {\it Example 2.} Consider the triangle singularities $X=D_{p,q,r}$
  with graph
  \[\Dfour{-r}{-p}{-1}{-q}
  \]
  It is known that they are minimally elliptic (i.~e., Gorenstein with
  $p_g(X)=1$, see \cite{Lauell})
  and have $\ep=1$, \cite{Wjac}. Hence $G=\C^*\times\bar{G}$.
  As $V(X)$ is spanned by the nowhere zero form $\omega$ we see that
  $X/H$ is minimally elliptic for $H\sub\bar{G}$ but $p_g(X/H)=0$
  (i.~e., $X/H$ is rational) if $H\not\sub\bar{G}$. Let us determine
  $X/\bar{G}$
  in case $\bar{G}\simeq\Aut\:\G$ is not trivial. First take $X=D_{p,q,q}$
  with $p\not=q$. Then $\bar{G}$ is generated by the involution $\s_3$ of
  Proposition 1. Hence $D_{p,q,q}/\bar{G}\simeq D_{2p,2,q}$ by an application
  of Lemmas 1 and 2.
  For $X=D_{p,p,p}$ we have $\bar{G}\simeq{\rm S}_3$. Let $H\le\bar{G}$ be
  the cyclic subgroup of order three. As $H$ acts freely on $\tilde{X}-E$
  the two fixed points lying on $E_0$ must be isolated. Lemma 4 below
  shows $D_{p,p,p}/H\simeq D_{3,3,p}$. Then
  $D_{p,p,p}/\bar{G}\simeq D_{3,3,p}/\s_3\simeq D_{2,3,2p}$.
\end{trivlist}
\begin{trivlist}
  \item[] \bf Lemma 4. \it Let $X=X_{b,1}$ with weighted graph
  \[\Aone{-b}
  \]
  Let $\s$ be an automorphism of $X$ of order three and consider its
  action on $\tilde{X}$ and on $E$. If $\s$ fixes $E$ pointwise then
  $X/\s$ has weighted graph
  \[\Aone{-3b}
  \]
  Otherwise there are the following possibilities:
  \[\Afour{}{}{-(b+3)/3}{-3}
  \qquad\mbox{or}\qquad
  \Aone{-b/3}
  \]
  if
  $b\equiv0$ {\rm mod} $3$,
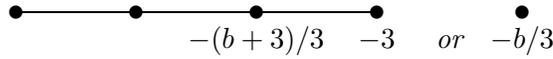
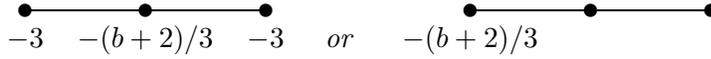
  \[\Athree{-3}{-(b+2)/3}{-3}
  \qquad\mbox{or}\qquad\qquad
  \Athree{-(b+2)/3}{}{}
  \]
  if
  $b\equiv1$ {\rm mod} $3$, and
  \[\Afive{}{}{-(b+4)/3}{}{}
  \qquad\mbox{or}\qquad\qquad
  \Atwo{-(b+1)/3}{-3}
  \]
  if
  $b\equiv2$ {\rm mod} $3$.
  The number of isolated fixed points lying on $E$ can be seen from the
  graphs.
  \item[] Proof. \rm
  Similar to the proof of Lemma 1.
  \hfill$\Box$
\end{trivlist}
\begin{trivlist}
  \item[] {\it Example 2 (continued).} Consider the fourteen triangle
  singularities which can be embedded as a hypersurface in $(\C^3,0)$,
  see \cite[Theorem 3.13 and section V]{Lauell}. In eleven cases only
  one of the
  three weights of the $\C^*$-action is odd so that the involution
  $\s_1\in\C^*$ is a reflection. By \cite[4.2]{MDiss} the quotient
  $X/\s_1$ will be smooth or an isolated hypersurface singularity, hence
  (as it is rational) a rational double point.
  For $D_{3,3,6}=Q_{12}$ and $D_{3,4,5}=S_{12}$ the $\C^*$-action has
  exactly two odd weights. The quotient $D_{3,3,6}/\s_1$ is a rational
  triple point with graph
  \[\Esix{-3}{}{}{}{}{}
  \]
  And $D_{3,4,5}/\s_1$ coincides with the icosahedral quotient singularity
  $I_{13}$ of graph
  \[\Esix{}{}{-3}{}{}{}
  \]
  This seems to be a quite interesting singularity: It plays an exceptional
  role in Manetti's \cite{Ma} study of smooth curves on rational surface
  singularities.
  Finally, for $X=D_{3,3,5}=Z_{13}$ all three weights are odd. Hence
  $X\to X/\s_1$ is unramified outside the singular point. In fact, it is
  the canonical Gorenstein cover \cite[section 4]{Wrat} of the rational
  quadruple point $X/\s_1$ with graph
  \[
  \unitlength0.8mm
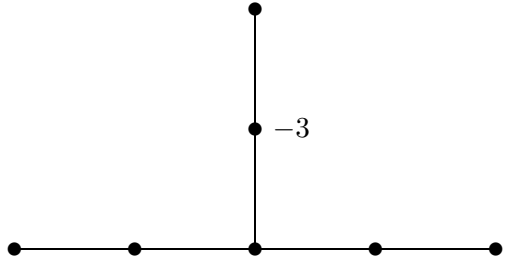
  \begin{picture}(80,49)
  \multiput(0,6)(20,0){5}{\circle*{2}}
  \multiput(40,26)(0,20){2}{\circle*{2}}
  \put(43,26){\makebox(0,0)[l]{$-3$}}
  \put(0,6){\line(1,0){80}}
  \put(40,6){\line(0,1){40}}
  \end{picture}
  \]
\end{trivlist}
\begin{trivlist}
  \item[] {\it Example 3.}
  Consider the fourteen exceptional families $f_a=f_0+a\cdot g$, $a\in\C$,
  of unimodal singularities, \cite[part II]{AGV}.
  Here the $f_0$ are weighted homogeneous and
  define exactly the fourteen triangle singularities $X_0$ which can be
  embedded as a hypersurface in $(\C^3,0)$.
  For $a\not= 0$ the $f_a$ are semi-weighted homogeneous and
  define a singularity $X_1$ whose analytic type is independent of $a$.
  The family $f_a$ is topologically trivial. Hence $X_0$ and
  $X_1$ have the same resolution graph $\G$, \cite[Theorem 2]{N}. We claim
  that the maximal reductive automorphism group of $X_1$ is
  $\Z_2\times\Aut\:\G$. For $i=0,1$ let $H_i^*$ be the group of
  $\phi\in\GL(3,\C)$ with $\phi f_i=c\cdot f_i$ for some $c\in\C^*$. By
  looking at the polynomial $f_0$ one sees
  that there is a subgroup $B\sub\GL(3,\C)$, isomorphic to $\Aut\:\G$,
  centralized by $\C^*$, intersecting $\C^*$ trivially and such that
  $\C^*\times B\sub H_0^*$. It follows from Proposition 2 in section 5
  below that $\C^*\times B=H_0^*$.
  Looking at $g$, which is in fact a monomial of (weighted) degree
  $\deg\:g=\deg\:f_0+2$, one sees that the subgroup $\Z_2\times B$ is
  contained in $H_1^*$. Now let $G$ be a maximal reductive automorphism
  group of $X_1$ acting on $(\C^3,0)$ by contact equivalences of $f_1$
  and containing $\Z_2\times B$.
  For each $\phi\in G$ there is a unit $u$, say with constant term $c$,
  such that
  \begin{equation}
  f_1\circ\phi=u\cdot f_1=c\cdot f_0+c\cdot g+
  \;\mbox{terms of degree}\; \ge d+3
  \end{equation}
  where $d=\deg\:f_0$. Here we have used that the weights $w_i$ are $\ge 3$.
  Let $\phi^i$ be the components of $\phi$ and write
  $\phi^i=\sum\phi_{ij}$ where $\deg\:\phi_{ij}=w_i+j$.
  It follows from $f_1\circ\phi=u\cdot f_1$ via
  \cite[Theorem 2.1]{GHP} that in the sum for $\phi^i$ only indices $j\ge0$
  occur. Writing $\phi_j=(\phi_{1j},\phi_{2j},\phi_{3j})$ we obtain
  \begin{equation}
  f_1\circ\phi=f_0\circ\phi_0+(\partial_xf_0\circ\phi_0)\cdot\phi_1
  +\;\mbox{terms of degree}\; \ge d+2
  \end{equation}
  and $f_0\circ\phi_0=c\cdot f_0$. In most of the fourteen cases (namely,
  if $\Aut\:\G=1$) the only monomial of degree $w_i$ is $x_i$. But also in
  the remaining cases one easily sees that $f_0\circ\phi_0=c\cdot f_0$
  forces $\phi_0$ to be linear. Therefore $\phi\mapsto\phi_0$ defines a
  homomorphism $G\to H_0^*=\C^*\times B$. It is not clear that $\phi_0$ is
  the linear part of $\phi$. But the usual linearization trick produces an
  analytic map germ $\psi$ with $\psi_0=1$ and
  $\phi_0\circ\psi=\psi\circ\phi$ for all $\phi\in G$. It follows from
  $\psi_0=1$ that the linear part of $\psi$ is triangular with trivial
  diagonal (if the weights are suitably ordered) and hence that $\psi$
  is invertible. Consequently $G$ is mapped isomorphically onto its image
  in $\C^*\times B$. It remains to show that this image is contained in
  $\Z_2\times B$. As $G$ contains $B$ it is enough to consider
  $\phi\in G$ with $\phi_0\in\C^*$, say
  $\phi_0(x)=t\cdot x=(t^{w_1}x_1,t^{w_2}x_2,t^{w_3}x_3)$.
  Comparing (1) and (2) once more we obtain
  \[0=(\partial_xf_0(t\cdot x))\cdot\phi_1
  =\sum_it^{d-w_i}\cdot\phi_{i1}\cdot\partial_{x_i}f_0.
  \]
  Since the partials of $f_0$ form a regular sequence in $\O_3$ we conclude
  that the $\phi_{i1}$ are contained in the Jacobian ideal $j(f_0)$.
  In each of the fourteen cases, the degree of every partial is larger than
  $w_i+1$ for all $i$.
  Hence $\phi_1=0$. Now we have a more precise
  version of (2):
  \[f_1\circ\phi=f_0(t\cdot x)+g(t\cdot x)
  +(\partial_xf_0(t\cdot x))\cdot\phi_2+
  \;\mbox{terms of degree}\;\ge d+3.
  \]
  We obtain
  $c\cdot g=t^{d+2}\cdot g+(\partial_xf_0(t\cdot x))\cdot\phi_2$. The
  second summand must vanish because otherwise $g\in j(f_0)$ which is
  not the case. Then $c\cdot g=t^{d+2}\cdot g$ and
  $c\cdot f_0=f_0(t\cdot x)=t^d\cdot f_0$ imply $t=\pm1$ proving the claim.
\end{trivlist}
\begin{trivlist}
  \item[] {\it Example 4.}
  Consider the quadrilateral singularities $X$ with graph
  \[
  \unitlength0.8mm
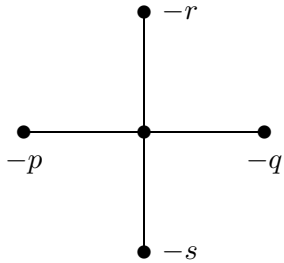
  \begin{picture}(40,43)
  \multiput(0,20)(20,0){3}{\circle*{2}}
  \multiput(20,0)(0,40){2}{\circle*{2}}
  \put(0,14){\makebox(0,0)[b]{$-p$}}
  \put(40,14){\makebox(0,0)[b]{$-q$}}
  \put(23,40){\makebox(0,0)[l]{$-r$}}
  \put(23,0){\makebox(0,0)[l]{$-s$}}
  \put(0,20){\line(1,0){40}}
  \put(20,0){\line(0,1){40}}
  \end{picture}
  \]
  Of course, the analytic type depends on the cross ratio of the four
  intersection points on the central curve $E_0\simeq\P_1$.
  The symmetric group ${\rm S}_4$ acts on
  $\C-\{0,1\}$ via the cross ratio. The subgroup $\Z_2\times\Z_2$
  consisting of $1$ and the three products of two disjoint transpositions
  acts trivially. The quotient map is
  \[j:\C-\{0,1\}\to\C:\l\mapsto
  \frac{4(\l^2-\l+1)^3}{27\l^2(\l-1)^2}.
  \]
  In particular, if the four intersection points on the central curve
  have cross ratio $\l$ (with respect to some numbering) then
  $j(X)=j(\l)$ is an invariant of the singularity (independent of the
  numbering). The generic orbit of ${\rm S}_4$ on $\C-\{0,1\}$ has
  cardinality six. There are two exceptional orbits: One of them,
  corresponding to $j=0$, consists of the two solutions of
  $\l^2-\l+1=0$, and the other, corresponding to $j=1$, consists of
  $-1$, $2$ and $1/2$. These well known facts are enough to determine,
  in each case, the group
  $A\simeq G/G_1$ of automorphisms
  of $\P_1$ which permute the intersection points while preserving the
  self intersection numbers. If $p=q=r=s$ then $A$ coincides with the
  alternating group ${\rm A}_4$ for $j(X)=0$. For $j(X)=1$ it is dihedral
  of order
  $8$ and for $j(X)\not=0,1$ it is $\Z_2\times\Z_2$.
  If $p=q=r\not=s$ then $A\simeq\Z_3$ for $j(X)=0$, $A\simeq\Z_2$
  generated by a transposition for
  $j(X)=1$, and $A=1$ for $j(X)\not=0,1$.
  For the remaining cases let $\l$ be the cross ratio of the four
  intersection points $p_1,\ldots,p_4$ (in this order) and suppose that
  the curves intersecting in $p_1$ and $p_2$ have equal self intersection
  number $p=q$ whereas $r$ and $s$ are different from $p$. There is an
  automorphism of $\P_1$ interchanging $p_1$ and $p_2$ while fixing $p_3$
  and $p_4$ if and only if $\l=-1$. If $r\not=s$ we see $A\simeq\Z_2$
  for $\l=-1$ and $A=1$ else.
  If $r=s$
  and $\l=-1$ then $A\simeq\Z_2\times\Z_2$ is generated by the two
  transpositions interchanging $p_1$ and $p_2$, respectively $p_3$ and $p_4$.
  Finally, if $r=s$ and $\l\not=-1$
  then $A\simeq\Z_2$ is generated by the product of those two
  transpositions.
  The quadrilateral singularities
  are minimally elliptic with $\ep=1$, \cite{Lauell} and \cite{Wjac}.
  Hence $G=G_1\times\bar{G}$ is a direct product, $X/H$ is
  rational for finite subgroups $H\le G$ with $H\not\sub\bar{G}$ and
  $X/H$ is minimally elliptic if $H\sub\bar{G}$.
  The reader is invited to determine the quotient $X/\bar{G}$. In each
  case, it is a triangle or a quadrilateral singularity.
\end{trivlist}

\section{The action on the homology of the link}
The link of a weighted homogeneous surface singularity $(X,0)$ is
a deformation retract of $X-0$ with $X$ denoting the corresponding affine
algebraic variety. The group $G$ acts on $X-0$.
Since the connected subgroup $G_1$ acts trivially on
integral homology there is an action of $G/G_1$ on
$H_1(X-0,\Z)$.
\begin{trivlist}
  \item[] \bf Theorem 6. \it
  Let $X$ be weighted homogeneous.
  Then the group $G/G_1$ acts faithfully on $H_1(X-0,\Z)$.
  \item[] Proof. \rm
  In the sequel all homology and cohomology modules are with integral
  coefficients. We need to recall how the homology of the link is expressed
  in terms of resolution data, see e.~g.\ \cite[section 4]{LW}.
  Let $s$ be the number of components of the exceptional
  divisor $E$ of the minimal good resolution $\tilde{X}$ and let $g$ be the
  genus of the central curve $E_0$. (We suppose that $X$ is not a cyclic
  quotient singularity. For those the proof is left to the reader.)
  Since $E$ is a deformation retract
  of $\tilde{X}$ we have $H_2(\tilde{X})\simeq H_2(E)\simeq \Z^s$ and
  $H_1(\tilde{X})\simeq H_1(E)\simeq H_1(E_0)\simeq \Z^{2g}$. Write
  $L=\tilde{X}-E\simeq X-0$. Then Lefschetz Duality gives
  $H_1(\tilde{X},L)\simeq H^3(E)=0$ and
  $H_2(\tilde{X},L)\simeq H^2(E)\simeq H_2(E)'$ with the prime denoting the
  dual $\Z$-module. Observe that all identifications are equivariant with
  respect to $G/G_1$.
  The exact homology sequence of the pair
  $(\tilde{X},L)$ comes down to
  \[H_2(E)\stackrel{j}{\to}H_2(E)'\to H_1(L)\to H_1(E_0)\to 0
  \]
  yielding the exact sequence
  \[0\to{\rm coker}\:j\to H_1(L)\to H_1(E_0)\to 0.
  \]
  Here $j$ denotes the adjoint of the intersection product:
  \[j(a):b\mapsto a\cdot b \quad\mbox{for}\quad a,b\in H_2(E).
  \]
  As the intersection matrix is negative definite the discriminant group
  ${\rm coker}\:j$ is torsion and hence equals the torsion subgroup $H_1(L)_t$
  of $H_1(L)$. Recall from \cite[V.3.1]{FK} that $\Aut\: E_0$ acts
  faithfully on $H_1(E_0)$ for $g\ge2$. This is not true for $g=1$ but then
  it is easily seen that the group $\Aut_{p_0}\, E_0$ of automorphisms
  fixing a point $p_0\in E_0$ acts faithfully on $H_1(E_0)$.
  Now consider the natural homomorphism $G/G_1\to\Aut\:\G$ and
  let $\phi\in G/G_1$ be trivial on $H_1(L)$. We claim that
  triviality on $H_1(L)_t$ forces $\phi$ to act trivially on $\G$.
  Accepting this for the moment, we are done in case
  $g=0$ since then $G/G_1\to\Aut\:\G$ is injective, see
  Remark (i) of section 2. For $g\ge1$ we conclude that $\phi$, viewed as
  an automorphism of the central curve $E_0$, has to fix the $r$ intersection
  points of $E_0$ with other components of $E$. As $\phi$ is trivial on the
  free part $H_1(E_0)$ of $H_1(L)$, too, we obtain $\phi=1$ in all cases
  except for simple elliptic $X$. Then,
  by Theorem 3 we have
  $G/G_1\simeq(\Z_b\times\Z_b)\rtimes\Aut_0\, E_0$
  with $-b=E_0\cdot E_0$. It is known \cite[pp.~282 - 283]{LW}
  that $\Z_b\times\Z_b$ is mapped isomorphically onto the group
  ${\rm Hom}(\Z^2,\Z_b)$ of automorphisms of $H_1(L)$ which act trivially
  on both $H_1(L)_t\simeq\Z_b$ and $H_1(E_0)\simeq\Z^2$. As $\Aut_0\, E_0$
  acts faithfully on $H_1(E_0)$ we conclude that $G/G_1$ acts faithfully
  on $H_1(L)$.
  \item[] Let us prove the claim. For $i=1,\ldots,r$ let $E_{ij}$,
  $j=1,\ldots,l_i$, be the curves on the $i$-th branch of the exceptional
  divisor (counted beginning at the centre), and $-b_{ij}=E_{ij}\cdot E_{ij}$
  as in section 3. Moreover, write $E_{00}=E_0$ and $-b=E_0\cdot E_0$.
  Let $\l_{ij}$ be the basis of $H_2(E)'$ dual to the basis $E_{ij}$ of
  $H_2(E)$. The action of $G/G_1$ on $H_2(E)$ is given by permutation
  of the curves: $\phi E_{ij}=E_{\phi(i,j)}$. Thus we have
  \[(\phi\l_{ij})(a)=a_{\phi(i,j)}\quad\mbox{for all}\quad
  a=\sum a_{ij}E_{ij}\in H_2(E).
  \]
  Triviality of $\phi$ on $H_1(L)_t={\rm coker}\: j$ means
  \[\phi\l-\l\in\im\: j\quad\mbox{for all}\quad\l\in H_2(E)'.
  \]
  Now assume that $\phi$ is not trivial on $\G$. Then one may assume that
  $\phi E_{11}=E_{21}$. There is $x=\sum x_{ij}E_{ij}\in H_2(E)$ such
  that $\phi\l_{11}-\l_{11}=j(x)$, i.~e., $x\cdot a=a_{21}-a_{11}$ for all
  $a\in H_2(E)$. This gives $x\cdot E_{11}=-1$, $x\cdot E_{21}=1$, and
  $x\cdot E_{ij}=0$ for all other $(i,j)$. By looking at the curves on the
  $i$-th branch one obtains $x_{i,l_i-1}=b_{i,l_i}\cdot x_{i,l_i}$, then
  $x_{i,l_i-2}=b_{i,l_i-1}\cdot x_{i,l_i-1}-x_{i,l_i}
  =(b_{i,l_i-1}-1/b_{i,l_i})\cdot x_{i,l_i-1}$ and so on up to
  \[x_{00}=\a_i/\b_i\cdot x_{i1}+\g_i
  \]
  where $\a_i/\b_i$ is the Hirzebruch-Jung continued fraction of the $i$-th
  branch and $\g_i=-1$, $1$ or $0$ according to $i=1$, $2$ or else.
  Finally, $x\cdot E_0=0$ implies
  \[0=-b\cdot x_{00}+\sum_{i=1}^r x_{i1}
  =(\sum_{i=1}^r \b_i/\a_i -b)\cdot x_{00}
  \]
  because $(\a_1,\b_1)=(\a_2,\b_2)$. The intersection matrix being negative
  definite we have $b>\sum_{i=1}^r\b_i/\a_i$ and $x_{00}=0$. Then
  $1=\a_1/\b_1\cdot x_{11}$ with $\a_1/\b_1>1$ and $x_{11}\in\Z$ yields
  a contradiction proving the claim.
  \hfill$\Box$
\end{trivlist}
\begin{trivlist}
  \item[] \it Remark. \rm
  The proof shows that $G/G_1$ even acts faithfully on the torsion
  subgroup of $H_1(X-0,\Z)$ if $r>2g+2$.
\end{trivlist}
Consider a weighted homogeneous normal surface singularity $X$ which can be
embedded as a hypersurface in $(\C^3,0)$ but which is not an
$A_k$-singularity. We may choose coordinates such that its maximal
reductive automorphism group $G$ acts linearly on $\C^3$ and such that
$\C^*=G_1\le G$ acts diagonally. Then $X$ is defined by a weighted
homogeneous polynomial $f$, say of degree $d$.
\begin{trivlist}
  \item[] \bf Proposition 2. \it
  In this situation $G$ equals the group $H^*$ of $\phi\in\GL(3,\C)$ with
  $\phi f=c\cdot f$ for some $c\in\C^*$. The subgroup $H$ of
  $\phi\in\GL(3,\C)$ with $\phi f=f$ is finite. The intersection
  $H\cap G_1=H\cap\C^*$ is the cyclic group $\Z_d$ of $d$-th roots of unity
  and
  $G/G_1\simeq H/\Z_d$.
  \item[] Proof. \rm
  First look at $H$. This is an algebraic subgroup of $\GL(3,\C)$. So one has
  to show that its Lie algebra $\bf h$ consisting of all derivations
  $D=\sum \l_i\partial_{x_i}$ with linear forms $\l_i$ such that
  $Df=0$ is reduced to $0$. For $f\in m^3$ this is shown by a standard
  argument which can be found at several places in the literature, e.~g.\
  in \cite{OS}. Now suppose that $f\notin m^3$. Since $f$ does not define
  an $A_k$-singularity we may assume $f=x_1^2+g$ with $g\in m^3$. The argument
  just mentioned shows that every $D\in\bf h$ is of form
  $D=x_1\sum_{i>1}a_i\partial_{x_i}$ with $a_i\in\C$. But then $Df=0$
  clearly implies $D=0$. Now turn to $H^*$. This is also an algebraic
  subgroup being the image under the projection
  $\GL(3,\C)\times\C^*\to\GL(3,\C)$ of the algebraic group of
  all pairs $(\phi,c)$ with $\phi f=c\cdot f$. Its Lie algebra consisting
  of all
  $D=\sum\l_i\partial_{x_i}$ such that $Df\in\C\cdot f$ is one dimensional
  spanned by the Euler derivation. As $\C^*\sub H^*$ this implies that
  $H^*$ is reductive. Because $G$ is a group of linear contact equivalences
  of $f$ and because $\C^*$ is central in $G$ one has $G\sub H^*$. But
  then $G=H^*$ by maximality. The remaining assertions are obvious.
  \hfill$\Box$
\end{trivlist}
\begin{trivlist}
  \item[] {\it Example 5.} Let $X$ be defined by $f=x_1^d+x_2^d+x_3^d$ with
  $d\ge 3$ and let $H,H^*$ be as above. Since the three weights of the
  $\C^*$-action are equal to $1$ any automorphism of $(\C^3,0)$
  commuting with $\C^*$ is linear in the given coordinates. This shows
  that $H^*$ is the maximal reductive automorphism group of $X$, i.\ e.,
  the coordinates are well chosen, compare Remark (vii) of section 2. We
  clearly have $\Z_d^3\rtimes{\rm S}_3\sub H$. To
  prove equality take $\phi\in\GL(3,\C)$ with $\phi f=f$. The ideal
  $((x_1x_2x_3)^{d-2})$ generated by the Hessian of $f$ is $\phi$-stable.
  Consequently $(x_1x_2x_3)$ is $\phi$-stable and, after permutation,
  the ideals $(x_1)$, $(x_2)$ and $(x_3)$ must be $\phi$-stable. This
  shows $H=\Z_d^3\rtimes{\rm S}_3$. We conclude that
  $G/G_1\simeq H/\Z_d$ has order $6d^2$. When $d$ is a multiple of $3$ we
  obtain examples where $G\to G/G_1$ does not admit a section. In fact,
  then the third root of unity $\zeta\in\C^*\le G$ is contained in the
  commutator subgroup of $G$, namely
  $\zeta=\phi\s\phi^{-1}\s^{-1}$ with
  $\phi(x_1,x_2,x_3)=(\zeta^2x_1,\zeta x_2,x_3)$ and
  $\s(x_1,x_2,x_3)=(x_2,x_3,x_1)$. Clearly this prevents $G$ from being a
  direct product of $\C^*$ and some subgroup.
  (Note that for $d=3$ we are discussing the simple elliptic singularity
  obtained by contracting the zero-section of a line bundle of degree
  $-3$ on the elliptic curve of $j$-invariant $0$.)
  On the contrary, if $3$ does not divide $d$ then $G$ is a direct product
  over $\C^*$. To prove this consider the normal subgroup
  $N=\Z_d^3\cap{\rm SL}(3,\C)$ of $\Z_d^3$. It has trivial intersection
  with the center $\Z_d$ of $H$ because $d$ is not a multiple of $3$.
  Then $|N|=d^2$ implies $\Z_d^3=\Z_d\times N$. As $N$ is
  ${\rm S}_3$-invariant we conclude $H=\Z_d\times B$ for
  $B=N\rtimes {\rm S}_3$ and then $G=\C^*\times B$.
\end{trivlist}
Returning to the general situation as described in Proposition 2,
choose an $H$-invariant Hermitean inner product on $\C^3$ and let
$\bar{B}_\ep$ be the corresponding closed ball of small radius $\ep$. One
has the Milnor fibration
\[f^{-1}(\bar{D}_\delta-0)\cap\bar{B}_\ep\to \bar{D}_\delta-0
\]
where $\bar{D}_\delta\sub\C$ is a small closed disc. Then clearly
the group $H$ of
right equivalences of $f$ acts on the Milnor fibre $F$. Observe that by
an equivariant version \cite[Lemma 4]{DEhr} of the Ehresmann Fibration
Theorem any two Milnor fibres are $H$-equivariantly diffeomorphic.
Moreover, let $M=H_2(F,\Z)$ equipped with the intersection form be the
Milnor lattice, $O(M)$ its group of isometries, and $W\le O(M)$ the
monodromy group.
\begin{trivlist}
  \item[] \bf Theorem 7. \it
  The homomorphism $H\to O(M)$ is injective and $H\cap W\simeq\Z_d$,
  generated by the monodromy operator.
  \item[] Proof. \rm
  Consider the Jacobian ideal $j(f)$ and $U=\O_3/j(f)$. Clearly $H$ acts
  on $U$. Wall \cite{Wsym}, see also \cite{OS}, has constructed an
  isomorphism $H_2(F,\C)\simeq U'\otimes \Lambda^3\C^3$ of $H$-modules.
  Let $\eta\in H$ be trivial on $M=H_2(F,\Z)$. As
  the basis element $1$ of $U$ is an eigenvector of eigenvalue $1$
  we conclude
  that $\det\:\eta=1$ and hence that $\eta$ must be trivial on $U$.
  Consequently $\eta x_i\equiv x_i$ mod $j(f)$ for the coordinate functions
  $x_i$. Because $X$ is not an $A_k$-singularity we may assume that
  $x_1$ is the only linear form contained in $j(f)$. Hence $\det\:\eta=1$
  implies
  \[\eta=\left(
    \begin{array}{ccc} 1 & * & * \\
                       0 & 1 & 0 \\
                       0 & 0 & 1
    \end{array}
    \right)
  \]
  Thus $\eta$, having finite order, must be trivial.
  \item[]
  As $X$ is defined by a
  weighted homogeneous polynomial of degree $d$ the $d$-th root of unity
  in $H$ induces a monodromy diffeomorphism $F\to F$. Therefore the
  monodromy operator has order $d$ and is contained in $H$ as acting on $M$.
  To show $H\cap W\sub \Z_d$ consider the exact homology sequence of the
  pair $(F,\partial F)$. By Lefschetz Duality it looks as follows:
  \[H_2(F,\Z)\stackrel{j}{\to}H_2(F,\Z)'\to H_1(\partial F,\Z)\to0
  \]
  where $j$ is the adjoint of the intersection form. The monodromy group
  $W$ is generated by Picard-Lefschetz transformations. These are
  reflections at hyperplanes orthogonal to the vanishing cycles. In
  particular, $W$ acts trivially on $M/M'$. Therefore $H\cap W$ acts
  trivially on $H_1(\partial F,\Z)$. Because the fibre bundle
  $f^{-1}(\bar{D}_\delta)\cap S_\ep\to \bar{D}_\delta$ is trivial and the link
  $L=X\cap S_\ep$ is a deformation retract of $X-0$ there are $H$-equivariant
  isomorphisms $H_1(\partial F,\Z)\simeq H_1(L,\Z)\simeq H_1(X-0,\Z)$.
  Then Theorem 6 implies that any $\eta\in H\cap W$ must be contained
  in $H\cap G_1=\Z_d$.
  \hfill$\Box$
\end{trivlist}

\vfill
Fachbereich Mathematik, Universit\"at Mainz, D 55099 Mainz, Germany \\
mueller@mat.mathematik.uni-mainz.de
\end{document}